\documentclass[12pt,sort&compress]{article}
\usepackage{amssymb}
\usepackage{amsfonts}
\usepackage{amsmath}
\usepackage{amssymb,amsmath}
\usepackage{amscd}
\def\a{\alpha}
\def\g{\gamma}

\def\Vir{\hbox{Vir}}
\def\cl{\centerline}

\def\vs{\vspace*}
\def\W{\mathcal{W}}

\def\Z{\mathbb{Z}}
\def\CC{\mathbb{C}}

\def\QED{\hfill$\Box$}

\def\t{\tilde}

\def\pa{\partial}
\def\la{\lambda}

\numberwithin{equation}{section}
\newtheorem{theo}{Theorem}[section]
\newtheorem{defi}[theo]{Definition}

\newtheorem{prop}[theo]{Proposition}

\newtheorem{rema}[theo]{Remark}

\makeatletter

\def\@biblabel#1{#1.~}

\makeatother

\begin{document}
\vs{10pt} \cl{\large {\bf Structures of $W(2,2)$ Lie conformal algebra}
\footnote{Corresponding author: Henan Wu (wuhenan@sxu.edu.cn).}} \vs{12pt}

\cl{ Lamei Yuan$^{\,\ddag}$, Henan Wu$^{\,\dag}$ }
 \cl{\small{ $^{\ddag}$ Academy of Fundamental and Interdisciplinary
 Sciences,}}\cl{\small{Harbin Institute of Technology, Harbin 150080, China}}
\cl{\small{$^{\dag}$School of Mathematical Sciences, Shanxi University, Taiyuan, 030006, China}}
\cl{\small E-mail:
lmyuan@hit.edu.cn, wuhenan@sxu.edu.cn
 }\vs{6pt}
{\small\parskip .005 truein \baselineskip 3pt \lineskip 3pt

\noindent{\bf Abstract:} The purpose of this paper is to study $W(2,2)$ Lie conformal algebra, which has a free $\mathbb{C}[\partial]$-basis $\{L, M\}$ such that
$[L_\lambda L]=(\partial+2\lambda)L$, $[L_\lambda
M]=(\partial+2\lambda)M$, $[M_\lambda M]=0$.  In this paper, we study conformal derivations, central extensions and conformal modules for this Lie conformal algebra. Also, we compute the cohomology of this Lie conformal algebra with coefficients in its modules. In particular, we determine its cohomology
with trivial coefficients both for the basic and reduced complexes.
 \vs{5pt}

\vs{5pt}

\noindent{\bf MR(2000) Subject Classification:}~17B10, 17B65, 17B68.


\vs{18pt}

\noindent{\bf 1. \
Introduction}\setcounter{section}{1}\setcounter{theo}{0}\setcounter{equation}{0}
\vs{6pt}

 A {\it Lie conformal algebra} is a $\CC[\partial ]$-module $\mathcal {R}$ equipped with a
$\lambda$-bracket $[\cdot_\lambda\cdot]$ which is a
$\CC$-bilinear map from $\mathcal {R}\otimes \mathcal {R}$ to $\CC[\lambda]\otimes \mathcal {R}$, such that the following axioms hold for all $a, b, c\in \mathcal {R}$:
\begin{eqnarray}
[\partial a_\lambda b]&=&-\lambda[a_\lambda b],\ \ [ a_\lambda \partial b]=(\partial+\lambda)[a_\lambda b] \ \ \mbox{(conformal\  sesquilinearity)},\label{Lc1}\\
{[a_\lambda b]} &=& -[b_{-\lambda-\partial}a] \ \ \mbox{(skew-symmetry)},\label{Lc2}\\
{[a_\lambda[b_\mu c]]}&=&[[a_\lambda b]_{\lambda+\mu
}c]+[b_\mu[a_\lambda c]]\ \ \mbox{(Jacobi \ identity)}\label{Lc3}.
\end{eqnarray}
In practice, the $\la$-brackets arise as generating functions for the singular part
of the operator product expansion in two-dimensional conformal field theory \cite{K1}. Lie conformal algebras are closely related to vertex algebras and infinite-dimensional Lie (super)algebras satisfying the locality property \cite{K3}. In the past few years, semisimple Lie conformal algebras have been intensively studied. In particular, the classification of all finite semisimple Lie conformal (super)algebras were given in \cite{DK,FK}. Finite irreducible conformal modules over the Virasoro, the current and the Neveu-Schwarz conformal algebras were classified in \cite{CK}. The cohomology theory was developed by  Bakalov, Kac, and Voronov in \cite{BKV}, where explicit computations of cohomologies of the Virasoro and the current conformal algebras were given. The aim of this paper is to study  structures including derivations, central extensions, conformal modules and cohomologies of a non-semisimple Lie conformal algebra associated to the $W$-algebra $W(2,2)$. The method here is based on the theory of the Virasoro conformal modules \cite{CK} and the techniques developed in \cite{DK, BKV}. Our discussions on the cohomology may be useful to do same thing for the other non-semisimple Lie conformal algebras, such as the Schr\"odinger-Virasoro, the loop Virasoro and the loop Heisenberg-Virasoro conformal algebras studied in \cite{FSW, SY,WCY}.

The Lie conformal algebra which we consider in this paper, denoted by $\mathcal{W}$, is a free $\mathbb{C}[\partial]$-module of  rank 2 generated by $L, M$, satisfying
\begin{eqnarray}\label{lamda-bracket}
[L_\lambda L]=(\partial+2\lambda)L,\ \ {[L_\lambda
M]}=(\partial+2\lambda)M,\ \
[M_\lambda M]=0.
\end{eqnarray}
The formal distribution Lie algebra corresponding to $\W$ is the centerless
$W$-algebra $W(2,2)$ introduced in \cite{ZD}. Besides, this Lie conformal algebra can be seen as a special case of a more general $W(a,b)$ Lie conformal algebra studied in \cite{XY}. Obviously, the Lie conformal algebra $\W$ contains the
Virasoro conformal algebra $\Vir$ as a subalgebra, namely,
\begin{eqnarray}\label{Vir}
\Vir=\CC[\partial]L,\ \ \ [L_\lambda L]=(\partial+2\lambda)L.
\end{eqnarray}
And it has a nontrivial
abelian conformal ideal generated by $M$. Thus it is not semisimple.

The rest of the paper is organized as follows. In Sect. 2, we study conformal derivations of the Lie conformal algebra $\mathcal
{W}$. It turns out that all the conformal derivations of $\mathcal
{W}$ are inner. In Sect. 3, we discuss
central extensions of $\W$ and show that $\W$ has a unique nontrivial
universal central extension. In Sect. 4,
we determine all free nontrivial conformal $\W$-modules of rank one. In Sect. 5, we compute cohomologies of $\W$ with coefficients in $\W$-modules $\CC$, $\CC_a$ and $M_{\Delta,\a}$, respectively. Consequently, we have the basic and reduced cohomologies for all $q\geq 0$ determined.

Throughout this paper, all vector spaces, linear maps and tensor products are over the complex field $\mathbb{C}$.  We use notations $\Z$ for the set of integers and $\Z_+$ for the set of nonnegative integers.
\vs{8pt}

 \noindent{\bf 2. Conformal derivation
}\setcounter{section}{2}\setcounter{theo}{0}\setcounter{equation}{0}\vs{8pt}

 Let $\mathcal {C}$
denote the ring $\mathbb{C}[\partial]$ of polynomials in the
indeterminate $\partial$.
\begin{defi}\rm Let $V$ and $W$ be two $\mathcal {C}$-modules.
A $\mathbb{C}$-linear map $\phi:V\rightarrow \mathcal
{C}[\lambda]\otimes_{\mathcal {C}}W$, denoted by $\phi_\la: V
\rightarrow W $, is called a  conformal linear map, if
$$\phi(\partial v)=(\partial+\lambda)(\phi v),\ \mbox{for }\ v\in V.$$
\end{defi}

 Denote by ${\rm Chom}(V,W)$ the space of conformal
linear maps from $\mathcal {C}$-modules $V$ to $W$. It can be made into an $\mathcal {C}$-module
via
\begin{eqnarray*}
(\pa\phi)_\la v=-\la\phi_\la v, \ \mbox{for} \ v\in V.
\end{eqnarray*}

\begin{defi}\rm Let $\mathcal {A}$ be a Lie conformal algebra.
A conformal linear map $d_\lambda:\mathcal {A}\rightarrow \mathcal
{A}$ is called a conformal derivation if
$$d_\lambda([a_\mu b])=[(d_\lambda a)_{\lambda+\mu}b]+[a_\mu(d_\lambda b)],\  \mbox{for all}\ a, b\in\mathcal {A}.$$
\end{defi}

 Denote by ${\rm CDer}(\mathcal {A})$ the space of all conformal derivations of $\mathcal {A}$. For any $a\in\mathcal {A}$, one can
define a linear map $({\rm ad}\,a)_\lambda:\mathcal
{A}\rightarrow \mathcal {A}$ by $({\rm ad}\,a)_\lambda
b=[a_\lambda b]$ for all $b\in\mathcal {A}.$ It is easy to
check that $({\rm ad}\,a)_\lambda$ is a conformal derivation of
$\mathcal {A}$. Any  conformal derivation of this kind is called an {\it inner derivation}. The space
of all inner derivations is denoted by ${\rm CInn}(\mathcal {A})$.

For the Lie conformal algebra $\mathcal
{W}$, we have the following result.
\begin{prop}\label{thm-2} Every conformal derivation of $\mathcal {W}$ is inner, namely,
${\rm CDer}(\mathcal {W})={\rm CInn}(\mathcal {W})$.
\end{prop}
\noindent {\it Proof.~}Let $d_\lambda$ be a conformal
derivation of $\mathcal {W}$. Then
\begin{eqnarray}\label{lm}
d_\lambda L=f_1(\lambda,\partial)L+f_2(\lambda,\partial)M,\ \
d_\lambda M=h_1(\lambda,\partial)L+h_2(\lambda,\partial)M,
\end{eqnarray}
where $f_i(\lambda,\partial)$ and $h_i(\lambda,\partial)$ for
$i=1,2$ are polynomials in $\CC[\lambda,\partial]$. Applying $d_\lambda$ to $[L_\mu L]=(\partial+2\mu) L$, we have
\begin{eqnarray*}\label{11_1}
d_\lambda[L_\mu L]&=&[(d_\lambda L)_{\lambda+\mu}L]+[L_\mu(d_\lambda
L)]\nonumber\\[4pt]
&=&(\partial+2\lambda+2\mu)f_1(\lambda,-\lambda-\mu)L+(\partial+2\lambda+2\mu)f_2(\lambda,-\lambda-\mu)M\nonumber\\
&&+(\partial+2\mu)f_1(\lambda,\partial+\mu)L+(\partial+2\mu)f_2(\lambda,\partial+\mu)M \\ &=&
d_\lambda((\partial+2\mu)
L)\\ &=&(\partial+\lambda+2\mu)\big(f_1(\lambda,\partial)L+f_2(\lambda,\partial)M\big).\label{22_2}
\end{eqnarray*}
Comparing the coefficients
of the similar terms gives
\begin{eqnarray}
(\partial+\lambda+2\mu)f_i(\lambda,\partial)-(\partial+2\mu)f_i(\lambda,\partial+\mu)=(\partial+2\lambda+2\mu)f_i(\lambda,-\lambda-\mu), \ \mbox{for} \  i=1,2.\label{f1}
\end{eqnarray}
Write $f_i(\lambda,\partial)=\sum_{j=0}^{n}a_{i,j}(\lambda)\partial^j$
with $a_{i,n}(\la)\neq 0$. Assume $n > 1$. Equating the coefficients of $\pa^n$ in (\ref{f1}) gives $(\lambda-n\mu)
a_{i,n}(\lambda)=0$. Thus $a_{i,n}(\la)=0$, a contradiction.
Therefore,
\begin{eqnarray}
f_i(\lambda,\partial)=a_{i,0}(\la)+a_{i,1}(\la)\pa, \ \mbox{for}\ i=1,2.
\end{eqnarray}
By replacing $d_\lambda$ by $d_\la-{\rm
ad\,}(a_{1,1}(-\pa)L)_\lambda-{\rm ad\,}(a_{2,1}(-\pa)M)_\lambda$, we can suppose $a_{1,1}(\la)=a_{2,1}(\la)=0.$ Then  plugging
$f_i(\lambda,\partial)=a_{i,0}(\la)$  into (\ref{f1}) gives $a_{i,0}(\la)=0$ for $i=1,2$. Thus $d_\la(L)=0$ by \eqref{lm}.
Fouthermore, applying $d_\lambda$ to $[L_\mu M]=(\partial+2\mu) M$, we have
\begin{eqnarray}\label{++}
(\partial+\lambda+2\mu)h_i(\lambda,\partial)=(\partial+2\mu)h_i(\lambda,\partial+\mu),\ \mbox{for}\  i=1,2.
\end{eqnarray}
Comparing the coefficients of highest degree of $\lambda$ in \eqref{++} gives $h_i(\lambda,\partial)=0$ for $i=1,2$. Hence
$d_\lambda (M)=0$ by \eqref{lm}. This concludes the proof.\QED
\begin{rema}\rm  Proposition \ref{thm-2} is equivalent to ${\rm H}^1 (\mathcal{W}, \mathcal{W})=0$.
\end{rema}
\vs{8pt}

\noindent{\bf 3. Central extension
}\setcounter{section}{3}\setcounter{theo}{0}\setcounter{equation}{0}
\vs{8pt}


An {\it
extension} of a Lie conformal algebra $\mathcal {A}$ by an abelian
Lie conformal algebra $\mathfrak a$ is a short exact sequence of Lie
conformal algebras
\begin{eqnarray*}
0\longrightarrow{\mathfrak a}\longrightarrow\hat{ \mathcal
{A}}\longrightarrow \mathcal{A}\longrightarrow 0.
\end{eqnarray*}
In this case $\hat{ \mathcal{A}}$ is also called an extension of
$\mathcal{A}$ by $\mathfrak a$. The extension is said to be {\it
central} if
\begin{eqnarray*}
{\mathfrak a}\subseteq Z(\hat{\mathcal{A} })=\{x\in \hat{\mathcal{A}
}\,\big|\, [x_\la y ]_{\hat{\mathcal{A} }}=0\ \mbox{for all}\ y\in
\hat {\mathcal{A}}\},\ \mbox{and}\ \pa{\mathfrak a}=0.
\end{eqnarray*}
Consider the central extension $\hat {\mathcal{A}}$ of
$\mathcal{A}$ by the trivial module $\CC$. This
means $\hat {\mathcal{A}}\cong \mathcal{A}\oplus \CC{\mathfrak c}$, and
\begin{eqnarray*}
[a_\la b]_{\hat {\mathcal{A}}}=[a_\la
b]_{\mathcal{A}}+f_\la(a,b){\mathfrak c}, \ \mbox{for} \
a,b\in\mathcal{A},\end{eqnarray*} where $f_\la:\mathcal{A}\times
\mathcal{A}\rightarrow \CC[\la]$ is a bilinear map. The axioms \eqref{Lc1}--\eqref{Lc3} imply the following properties of the $2$-cocycle $f_\la(a,b)$:
\begin{eqnarray}
f_\la(a, b)&=&-f_{-\la-\pa}(b, a),\ \label{2-cocycle1}\\
f_\la(\pa a, b)&=&-\la f_\la(a, b)=- f_\la(a,\pa b),\label{2-cocycle2}\\
f_{\la+\mu}([a_\la b],c)&=&f_\la(a,[b_\mu c])-f_\mu(b,[a_\la c]),
\label{2-cocycle3}
\end{eqnarray}
for all $a, b, c\in \mathcal{A}$. For any linear function $f:\mathcal{A}\rightarrow \CC$, the map
\begin{eqnarray}
\psi_{f}(a,b)=f([a_\la b]),\ \mbox{for}\ a,b\in\mathcal{A},
\end{eqnarray}
defines a trivial $2$-cocycle. Let $a'_{\la}(a,b)=a_{\la}(a,b)+\psi_{f}(a,b)$. The equivalent
$2$-cocycles $a'_{\la}(a,b)$ and $a_{\la}(a,b)$ define isomorphic extensions.

In the following we determine the central extension $\widehat {\W}$ of $\W$ by $\CC{\mathfrak c}$, i.e., $\widehat {\W}=\W\oplus \CC{\mathfrak c}$, and the relations in (\ref{lamda-bracket})
are replaced by
\begin{eqnarray}\label{e-be}
[L_\la L]=(\pa+2\la)L+a_\la(L,L){\mathfrak c},\ \ {[L_\la
M]}=(\pa+2\la)M+a_\la(L,M){\mathfrak c},\ {[M_\la
M]}=a_\la(M,M){\mathfrak c},
\end{eqnarray}
and the others can be obtained by skew-symmetry. Applying the Jacobi identity
for $(L, L, L)$, we have
\begin{eqnarray}
(\la+2\mu)a_\la(L,L)-(\mu+2\la)a_\mu(L,L)=(\la-\mu)a_{\la+\mu}(L,L).
\label{cl1}
\end{eqnarray}
Write $a_\la(L,L)=\sum_{i=0}^{i=n}a_i\la^i\in\CC[\la]$ with $a_n\neq
0$. Then, assuming $n > 1$ and equating the coefficients of $\la^n$
in (\ref{cl1}), we get $2\mu a_n=(n-1)\mu a_n $ and thus $n=3$. Then
$$a_\la(L,L)=a_0+a_1\la+a_2\la^2+a_3\la^3.$$ Plugging this in
(\ref{cl1}) and comparing the similar terms, we obtain
$a_0=a_2=0$. Thus
\begin{eqnarray}\label{a1===}
a_\la(L,L)=a_1\la+a_3\la^3.\end{eqnarray}
 To compute $a_\la(L,M)$, we apply the Jacobi identity for $(L,L,M)$ and
obtain
\begin{eqnarray*}
(\la+2\mu)a_\la(L,M)-(\mu+2\la)a_\mu(L,M)=(\la-\mu)a_{\la+\mu}(L,M).
\end{eqnarray*}
By doing similar discussions as those in the process of computing
$a_\la(L,L)$, we have
\begin{eqnarray}\label{b1===}
a_\la(L,M)=b_1\la+b_3\la^3, \ \mbox{for some} \ b_1, b_3\in\CC.
\end{eqnarray} Finally, applying the
Jacobi identity for $(L,M,M)$ yields
$(\la-\mu)a_{\la+\mu}(M,M)=-(2\la+\mu)a_\mu(M,M)$, which implies
\begin{eqnarray}\label{h1===}
 a_\la(M,M)=0.
\end{eqnarray}

From the discussions above, we obtain the following results.
\begin{theo}\label{thm-4}
\begin{itemize}\parskip-3pt
\item[\rm(1)]
For any $a, b\in\CC$ with
$(a,b)\ne(0,0)$, there exists a unique nontrivial universal central
extension of the Lie conformal algebra $\W$ by $\CC{\mathfrak c}$, such that
\begin{eqnarray}\label{*****}
[L_\la L]=(\pa+2\la)L+a \la^3{\mathfrak c},\ {[L_\la
M]}=(\pa+2\la)M+b\la^3{\mathfrak c},\ [M_\lambda M]=0.
\end{eqnarray}
\item[\rm(2)]
There exists a unique nontrivial universal central extension of $\W$
by $\CC{\mathfrak c}\oplus\CC{\mathfrak c}'$, satisfying
\begin{eqnarray}
[L_\la L]=(\pa+2\la)L+\la^3{\mathfrak c},\ {[L_\la
M]}=(\pa+2\la)M+\la^3{\mathfrak c'},\ [M_\lambda M]=0.
\end{eqnarray}
\end{itemize}
\end{theo}
\noindent{\it Proof.~}(1) By \eqref{a1===}--\eqref{h1===}, replacing $L, M$ respectively by
$L-\frac12a_1{\mathfrak c},\,M-\frac12b_1{\mathfrak c}$ and noticing
that $\partial{\mathfrak c}=0$, we can suppose $a_1=b_1=0$. This shows
\eqref{*****}. The universality of the extension
follows from \cite{MP} and the fact that $\W$ is perfect, namely, $[\W_{\la}\W]=\W[\lambda]$.

(2) This follows from the proof of (1).
\QED

\begin{rema}\rm  Theorem \ref{thm-4}(2) implies that ${\rm dim\, H}^2 (\mathcal{W}, \CC)=2$.
\end{rema}
\noindent{\bf 4. Conformal module
}\setcounter{section}{4}\setcounter{theo}{0}\setcounter{equation}{0} \vskip8pt

Let us first recall the notion of a conformal module given in \cite{CK}.
\begin{defi}\rm A (conformal) module $V$ over a Lie conformal algebra $\mathcal {A}$
is a $\mathbb{C}[\partial]$-module endowed with a bilinear map $\mathcal {A}\otimes V\rightarrow
V[\lambda]$, $a \otimes v\mapsto a_\lambda v$ satisfying the following axioms for $a,b\in\mathcal {A}$, $v\in V$:
\begin{eqnarray*}
&&a_\lambda(b_\mu v)-b_\mu(a_\lambda v)=[a_\lambda b]_{\lambda+\mu}v,\\
&&(\partial a)_\lambda v=-\lambda a_\lambda v,\ \  a_\lambda(\partial
v)=(\partial+\lambda)a_\lambda v.
\end{eqnarray*}
An $\mathcal{A}$-module $V$ is called {\it finite} if it is finitely generated over $\mathbb{C}[\partial]$.
\end{defi}

The vector space $\CC$ can be seen as a module (called
the {\it trivial module}) over any conformal algebra $\mathcal {A}$
with both the action of $\pa$ and the action of $\mathcal {A}$ being
zero. For a fixed nonzero complex constant $a$, there
is a natural $\CC[\pa]$-module $\CC_a$, which is the one-dimensional
vector space $\CC$ such that $\pa v=a v$ for $v\in\CC_a$. Then
$\CC_a$ becomes an $\mathcal{A}$-module, where
$\mathcal{A}$ acts as zero.

 For the Virasoro conformal algebra $\Vir$, it is known from
\cite{DK} that all the free nontrivial $\Vir$-modules of rank one
over $\mathbb{C}[\partial]$ are the following ones $(\Delta,
\alpha\in \mathbb{C})$:
\begin{eqnarray}
M_{\Delta,\alpha}=\mathbb{C}[\partial]v,\ \ L_\lambda
v=(\partial+\alpha+\Delta \lambda)v.
\end{eqnarray}
The module $M_{\Delta,\alpha}$ is irreducible if and only if
$\Delta\neq 0$. The module $M_{0,\alpha}$ contains a unique
nontrivial submodule $(\partial +\alpha)M_{0,\alpha}$ isomorphic to
$M_{1,\alpha}.$  Moreover, the modules $M_{\Delta,\alpha}$ with
$\Delta\neq 0$ exhaust all finite irreducible nontrivial
$\Vir$-modules.

The following result describes the free nontrivial $\W$-modules of rank one.
Similar result for the more general $W(a,b)$ Lie conformal algebra was given in \cite{XY}. We aim to consider it in details in the $W(2,2)$ case.
\begin{prop} \label{thm-3} All free nontrivial $\W$-modules of rank one over
$\mathbb{C}[\pa]$ are the following ones:
\begin{eqnarray*}
M_{\Delta,\alpha}=\mathbb{C}[\partial]v,\ L_\lambda
v=(\partial+\alpha+\Delta \lambda)v, \ M_\lambda v=0, \ \mbox{for
some}\ \Delta,\a\in\CC.
\end{eqnarray*}
\end{prop}
\noindent{\it Proof.~} Suppose that $L_\lambda v=f(\partial,\lambda)v,\ M_\lambda v=g(\partial,\lambda)v,$ where
$f(\partial,\lambda)$, $g(\partial,\lambda)\in\mathbb{C}[\lambda,\partial]$. By the result of
$\Vir$-modules, we have
\begin{eqnarray*}
f(\partial,\lambda)=\partial+\alpha+\Delta \lambda, \ \ \mbox{for
some} \ \alpha,\Delta\in\mathbb{C}.
\end{eqnarray*}
On the other hand, it follows from $M_\lambda(M_\mu v)=M_\mu(M_\lambda v)$ that
$$g(\partial,\lambda)g(\pa+\la,\mu)=g(\pa,\mu)g(\pa+\mu,\la).$$ This
implies ${\rm deg}_\la g(\pa,\la)+{\rm deg}_\pa g(\pa,\la)={\rm
deg}_\la g(\pa,\la)$, where the notation ${\rm deg}_\la g(\pa,\la)$
stands for the highest degree of $\la$ in $g(\pa,\la)$. Thus ${\rm
deg}_\pa g(\pa,\la)=0$ and so $g(\pa,\la)=g(\la)$ for
some $g(\la)\in\CC[\la]$. Finally, $[L_\la M]_{\la+\mu}v=(\la-\mu)
M_{\la+\mu}v$ gives $(\la-\mu) g(\la+\mu)=-\mu g(\mu),$ which
yields $g(\la)=0$. This proves the result. \QED

\vs{8pt} \noindent{\bf 5. Cohomology
}\setcounter{section}{5}\setcounter{theo}{0}\setcounter{equation}{0}
\vs{8pt}

For completeness, we recall the following definition from \cite{BKV}:
\begin{defi}\label{definition}\rm An $n$-cochain ($n\in\Z_+$) of a Lie conformal algebra $\mathcal{A}$ with coefficients in an
$\mathcal{A}$-module $V$ is a $\CC$-linear map\vs{-5pt}
\begin{eqnarray*}
\gamma:\mathcal{A}^{\otimes n}\rightarrow V[\la_1,\cdots,\la_n],\ \
\ a_1\otimes\cdots \otimes a_n \mapsto
\g_{\la_1,\cdots,\la_n}(a_1,\cdots,a_n)
\end{eqnarray*}
satisfying the following conditions:\begin{itemize}\parskip-3pt
\item[\rm(1)] $\g_{\la_1,\cdots,\la_n}(a_1,\cdots,\pa a_i,\cdots,
a_n)=-\la_i\g_{\la_1,\cdots,\la_n}(a_1,\cdots, a_n)$ \ (conformal antilinearity),
\item[\rm (2)] $\g$ is skew-symmetric with respect to simultaneous permutations
of $a_i$'s and $\la_i$'s \ (skew-symmetry).
\end{itemize}
\end{defi}

As usual, let $\mathcal{A}^{\otimes 0}= \CC$, so that a $0$-cochain
is an element of $V$. Denote by  ${\t C}^n(\mathcal {A},V)$ the set
of all $n$-cochains. The differential $d$ of an $n$-cochain $\g$ is
defined as follows:
\begin{eqnarray}\label{ddd}
&&(d\g)_{\la_1,\cdots,\la_{n+1}}(a_1,\cdots,a_{n+1})\nonumber\\&&\ \ \ =\mbox{$\sum\limits_{i=1}^{n+1}$}(-1)^{i+1}a_{i_{\la_i}}\g_{\la_1,\cdots,\hat{\la_i},\cdots,\la_{n+1}}(a_1,\cdots,\hat{a_i},\cdots,a_{n+1})\nonumber\\
&&\ \ \ \ \ \ \ +\mbox{$\sum\limits_{i,j=1\atop i<j}^{n+1}$}(-1)^{i+j}\g_{\la_i+\la_j,\la_1,\cdots,\hat{\la_i},\cdots,\hat{\la_j},\cdots,\la_{n+1}}([a_{i_{\la_i}}a_j],a_1,\cdots,\hat{a_i},\cdots,\hat{a_j},\cdots,a_{n+1}),
\end{eqnarray}
where $\g$ is linearly extended over the polynomials in $\la_i$. In
particular, if $\g\in V$ is a $0$-cochain, then
$(d\g)_\la(a)=a_\la\g$.

It is known from \cite{BKV} that the operator $d$ preserves the space of cochains and $d^2=0$. Thus the cochains of a Lie conformal algebra $\mathcal{A}$ with coefficients in its module $V$ form a complex, which is denoted by
\begin{eqnarray}
\t C^\bullet(\mathcal{A},V)=\mbox{$\bigoplus\limits_{
n\in\Z_+}$}\t C^n(\mathcal{A},V),
\end{eqnarray}
and called the {\it basic complex}. Moreover, define a (left) $\CC[\pa]$-module structure on $\t
C^\bullet(\mathcal{A},V)$ by
\begin{eqnarray*}
(\pa\g)_{\la_1,\cdots,\la_n}(a_1,\cdots,
a_n)=(\pa_V+\mbox{$\sum\limits_{i=1}^n$}\la_i)\g_{\la_1,\cdots,\la_n}(a_1,\cdots,
a_n),
\end{eqnarray*}
where $\pa_V$ denotes the action of $\pa$ on $V$. Then $d\pa=\pa d$
and thus $\pa \t C^\bullet(\mathcal{A},V)\subset \t
C^\bullet(\mathcal{A},V)$ forms a subcomplex. The quotient
complex
\begin{eqnarray*}
C^\bullet(\mathcal{A},V)=\t C^\bullet(\mathcal{A},V)/\pa \t
C^\bullet(\mathcal{A},V)= \mbox{$\bigoplus\limits_{n\in\Z_+}$}
C^n(\mathcal{A},V)
\end{eqnarray*}
is called the {\it reduced complex}.
\begin{defi}\label{def11}\rm The basic cohomology ${\rm \t H}^\bullet (\mathcal{A},V)$ of a Lie conformal algebra $\mathcal{A}$ with coefficients
 in an $\mathcal{A}$-module $V$ is the
cohomology of the basic complex $\t C^\bullet(\mathcal{A},V)$ and
the (reduced) cohomology ${\rm \,H}^\bullet (\mathcal{A},V)$ is the
cohomology of
the reduced complex $C^\bullet(\mathcal{A},V)$.
\end{defi}

For a $q$-cochain $\gamma\in{\t C}^q(\mathcal {A},V)$, we call
$\gamma$ a {\it $q$-cocycle} if $d(\gamma)=0$; a {\it $q$-coboundary} if there exists a $(q-1)$-cochain $\phi\in\t
C^{q-1}(\mathcal{A},V)$ such that $\gamma=d(\phi)$. Two cochains
$\gamma_1$ and $\gamma_2$ are called {\it equivalent} if $\gamma_1-\gamma_2$ is a
coboundary. Denote by $\t D^q(\mathcal{A},V)$ and $\t B^q(\mathcal{A},V)$ the spaces
of $q$-cocycles and $q$-boundaries, respectively. By Definition \ref{def11},
\begin{eqnarray*}
{\rm \t H}^q(\mathcal{A},V)=\t D^q(\mathcal{A},V)/\t B^q(\mathcal{A},V)=\{\mbox{equivalent classes of
$q$-cocycles}\}.
\end{eqnarray*}

The main results of this section are the following.
\begin{theo}\label{thm-5}  For the Lie conformal algebra $\mathcal{W}$, the
following statements hold.\begin{itemize}\parskip-3pt
\item[\rm(1)] For the trivial module $\CC$,
\begin{eqnarray}\label{main1}
{\rm dim\,\t H}^q(\mathcal{W},\CC)=\left\{
\begin{array}{ll}
1 &{\mbox if}\ q=0,4,5,6,\\
2 &{\mbox if}\ q=3,\\
0 &{\mbox otherwise},
\end{array}
\right.
\end{eqnarray}
and
\begin{eqnarray}\label{main2}
{\rm dim\, H}^q(\mathcal{W},\CC)=\left\{
\begin{array}{ll}
1 &{\mbox if}\ q=0, 6,\\
2 &{\mbox if}\ q=2, 4, 5\\
3 &{\mbox if}\ q=3,\\
0 &{\mbox otherwise}.
\end{array}
\right.
\end{eqnarray}
\item[\rm(2)] If $a\neq 0$, then  ${\rm dim\, H}^q(\W,\CC_a)=0$, for $q\geq0$.
\item[\rm(3)] If $\a\neq 0$, then
$ {\rm dim\, H}^q(\mathcal{W},M_{\Delta,\a})=0,$ for $q\geq0$.
\end{itemize}
\end{theo}
\noindent{\it Proof.~} (1) For any $\gamma\in \t
C^0(\mathcal{W},\CC)=\CC$, we have $(d\gamma)_\la (X)=X_\la \gamma =0$ for
$X\in \mathcal{W}$. This means $\t D^0(\mathcal{W},\CC)=\CC$ and $\t B^0(\mathcal{W},\CC)=0$. Thus ${\rm \t
H}^0(\mathcal{W},\CC)=\CC$ and $ {\rm H}^0(\mathcal{W},\CC)=\CC$
since $\pa\CC=0$.

Let $\gamma\in\t C^1(\mathcal{W},\CC)$ be such that
$d\gamma\in\pa\t C^2(\mathcal{W},\CC)$, namely, there is $\phi\in\t
C^2(\mathcal{W},\CC)$ such that
\begin{eqnarray}\label{5--1}
\g_{\la_1+\la_2}([X_{\la_1} Y])=-(d\g)_{\la_1,\la_2}(X,Y)=-(\pa
\phi)_{\la_1,\la_2}(X,Y)=-(\la_1+\la_2)\phi_{\la_1,\la_2}(X,Y),
\end{eqnarray}
for $X,Y\in\{L,M\}$. By \eqref{5--1} and \eqref{lamda-bracket},
\begin{eqnarray}
(\la_1-\la_2)\g_{\la_1+\la_2}(X)=-(\la_1+\la_2)\phi_{\la_1,\la_2}(L,X),\ \ X\in\{L,M\}.\label{5*}
\end{eqnarray}
Letting $\la=\la_1+\la_2$ in \eqref{5*} gives
\begin{eqnarray}
(\la-2\la_2)\g_\la(X)=-\la\phi_{\la_1,\la_2}(L,X),\ \ X\in\{L,M\},\label{7-1}
\end{eqnarray}
which implies that $\g_\la(X)$ is divisible by $\la$. Define
\begin{eqnarray*}
\g'_{\la}(X)=\la^{-1}\g_{\la}(X),\ \ X\in\{L,M\}.
\end{eqnarray*}
Clearly, $\g'\in\t C^1(\W,\CC)$ and $\g=\pa \g'\in\pa\t
C^1(\W,\CC)$. Thus ${\rm H}^1(\W,\CC)=0$. If $\g$ is
a 1-cocycle (this means $\phi=0$), then \eqref{7-1} gives
$\g=0$. Hence, ${\rm \t H}^1(\W,\CC)=0$.

Let $\psi$ be a 2-cocycle. For $X\in\W$, we have
\begin{eqnarray*}
0&=&(d\psi)_{\la_1,\la_2,\la_3}(X,L,L)\\
&=&-(\la_1-\la_2)\psi_{\la_1+\la_2,\la_3}(X,L)+(\la_1-\la_3)\psi_{\la_1+\la_3,\la_2}(X,L)-(\la_2-\la_3)\psi_{\la_2+\la_3,\la_1}(L,X).
\end{eqnarray*}
Letting $\la_3=0$ and $\la_1+\la_2=\la$ gives $(\la-2\la_2)\psi_{\la,0}(X,L)=\la \psi_{\la_1,\la_2}(X,L)$. Hence, $\psi_{\la,0}(X,L)$ is divisible by $\la$. Define
a $1$-cochain $f$ by
\begin{eqnarray}\label{f}
f_{\la_1}(L)=\la_1^{-1}\psi_{\la_1,\la}(L,L)|_{\la=0},\ \
f_{\la_1}(M)=\la_1^{-1}\psi_{\la_1,\la}(M,L)|_{\la=0}.
\end{eqnarray}
Set $\g=\psi+d f$, which is also a 2-cocycle. By \eqref{f},
\begin{eqnarray}
\g_{\la_1,\la}(L,L)|_{\la=0}&=&\psi_{\la_1,\la}(L,L)|_{\la=0}-\la_1f_{\la_1}(L)=0,\label{5+}\\
\g_{\la_1,\la}(M,L)|_{\la=0}&=&\psi_{\la_1,\la}(M,L)|_{\la=0}-\la_1f_{\la_1}(M)=0.\label{5++}
\end{eqnarray}
By \eqref{5+}, we have
\begin{eqnarray*}
0&=&(d\g)_{\la_1,\la_2,\la}(L,L,L)|_{\la=0}\\
 &=&-\g_{\la_1+\la_2,\la}([L_{\la_1}L],L)|_{\la=0}+\g_{\la_1+\la,\la_2}([L_{\la_1
}L],L)|_{\la=0}-\g_{\la_2+\la,\la_1}([L_{\la_2} L],L)|_{\la=0}\\
&=&\la_1\g_{\la_1,\la_2}(L,L)-\la_2\g_{\la_2,\la_1}(L,L)\\
&=&(\la_1+\la_2)\g_{\la_1,\la_2}(L,L).
\end{eqnarray*}
Thus $\g_{\la_1,\la_2}(L,L)=0$. Similarly, by \eqref{5++},
\begin{eqnarray*}
0=(d\g)_{\la_1,\la_2,\la}(L,M,L)|_{\la=0}=(\la_1+\la_2)\g_{\la_1,\la_2}(L,M),
\end{eqnarray*}
which gives $\g_{\la_1,\la_2}(L,M)$=0 and so $\g_{\la_1,\la_2}(M,L)$=0.
Finally,
\begin{eqnarray}
0=(d\g)_{\la_1,\la_2,\la}(L,M,M)|_{\la=0}
=-(\la_1-\la_2)\g_{\la_1+\la_2,0}(M,M)+\la_1\g_{\la_1,\la_2}(M,M).\label{7-2}
\end{eqnarray}
Setting $\la_1=0$ in \eqref{7-2} gives
$\g_{\la_2,0}(M,M)=0$ and thus
$\g_{\la_1,\la_2}(M,M)=0$. This shows $\g=0$. Hence ${\rm \t
H}^2(\W,\CC)=0$. According to Theorem \ref{thm-4}(2), ${\rm dim\, H}^2(\W,\CC)=2$.

To determine higher dimensional cohomologies (for $q\geq 3$), we define an operator $\tau:\t C^q(\W,\CC)\rightarrow
\t C^{q-1}(\W,\CC)$ by
\begin{eqnarray}\label{7++}
(\tau
\g)_{\la_1,\cdots,\la_{q-1}}(X_1,\cdots,X_{q-1})=(-1)^{q-1}\frac{\pa}{\pa\la}\g_{\la_1,\cdots,\la_{q-1},\la}(X_1,\cdots,X_{q-1},L)|_{\la=0},
\end{eqnarray}
for $X_1,\cdots,X_{q-1}\in\{L,M\}$. By \eqref{ddd}, \eqref{7++} and skew-symmetry of $\g$,
\begin{eqnarray}\label{7++1}
&&((d\tau+\tau d)
\g)_{\la_1,\cdots,\la_{q}}(X_1,\cdots,X_q)\nonumber\\
&&\ \ \ \ =(-1)^q\frac{\pa}{\pa\la}\mbox{$\sum\limits_{i=1}^q$}(-1)^{i+q+1}\g_{\la_i+\la,\la_1,\cdots,\hat{\la_i},\cdots,\la_{q}}([X_{i\,{\la_i}} L],X_1,\cdots,\hat{X_{i}},\cdots, X_q)|_{\la=0}\nonumber\\
&&\ \ \ \ =\frac{\pa}{\pa\la}\mbox{$\sum\limits_{i=1}^q$}\g_{\la_1,\cdots,\la_{i-1},\la_i+\la,\la_{i+1},\cdots,\la_{q}}
(X_1,\cdots,X_{i-1},[X_{i\,{\la_i}} L],X_{i+1},\cdots,X_q)|_{\la=0}.
\end{eqnarray}
By the fact that $[X_{i\,{\la_i}} L]=(\pa+2\la_i)X_i$ and conformal antilinearity of $\g$, $[X_{i\,{\la_i}} L]$ can be replaced by
$(\la_i-\la)X_i$ in \eqref{7++1}.
Thus, equality \eqref{7++1} can be rewritten as
\begin{eqnarray}\label{7++2}
&&((d\tau+\tau d)
\g)_{\la_1,\cdots,\la_{q}}(X_1,\cdots,X_q)\nonumber\\
&&\ \ \ \ =\frac{\pa}{\pa\la}\mbox{$\sum\limits_{i=1}^q$}(\la_i-\la)\g_{\la_1,\cdots,\la_{i-1},\la_i+\la,\la_{i+1},\cdots,\la_{q}}
(X_1,\cdots,X_{i-1},X_i,X_{i+1},\cdots,X_q)|_{\la=0}\nonumber\\
&&\ \ \ \ =({\rm deg\,} \g-q)\g_{\la_1,\cdots,\la_{q}}(X_1,\cdots,X_q),
\end{eqnarray}
where ${\rm deg\,}\g$ is the total degree of $\g$ in $\la_1,\cdots,\la_q$. As it was explained in \cite{BKV}, only those homogeneous cochains whose degree as a polynomial is equal to their degree as a cochain contribute to the cohomology of $\t C^\bullet(\mathcal{W},\CC)$. Without loss of generality, we always assume that the first $k$ variables are $L$ and the last $q-k$ variables are $M$ in $\g_{\la_1,\cdots,\la_{q}}(X_1,\cdots,X_q)$, so that $\g_{\la_1,\cdots,\la_{q}}(X_1,\cdots,X_q)$ as a polynomial in $\la_1,\cdots,\la_q$ is skew-symmetric in $\la_1,\cdots,\la_k$ and also skew-symmetric in $\la_{k+1},\cdots,\la_q$. Therefore, it is divisible by $$\mbox{$\prod\limits_{1\leq i< j\leq k}$}(\la_i-\la_j)\times\mbox{$\prod\limits_{k+1\leq i< j\leq q}$}(\la_i-\la_j),$$ whose polynomial degree is $k(k-1)/2+(q-k)(q-k-1)/2$. Consider the quadratic inequality $k(k-1)/2+(q-k)(q-k-1)/2\leq q,$ whose discriminant is $-4k^2+12k+9$. Since $-4k^2+12k+9\geq 0$ has $k=0, 1, 2$ and $3$ as the only integral solutions, we have
\begin{eqnarray}\label{qqq}
q=\left\{
\begin{array}{ll}
0, 1, 2, 3, & if \ k=0,\\
1, 2, 3, 4, & if \ k=1,\\
2, 3, 4, 5, & if \ k=2,\\
3, 4, 5, 6, & if \ k=3.
\end{array}
\right.
\end{eqnarray}
Thus ${\rm \t H}^q(\mathcal{W},\CC)=0$ for $q\geq 7$. It remains to compute ${\rm \t H}^q(\W,\CC)$ for $q=3,4,5,6$.

 For $q=3$, we need to consider four cases for $k$, i.e., $k=0, 1,2,3$.  Let $\g\in\t D^3(\W,\CC)$ be a 3-cocycle. A direct computation shows that
\begin{eqnarray*}\label{q=3-1}
0=(d\gamma)_{\la_1,\la_2,\la_3,\la}(M,M,M,L)|_{\la=0}
=-(\la_1+\la_2+\la_3)\gamma_{\la_1,\la_2,\la_3}(M,M,M).
\end{eqnarray*}
This gives $\gamma_{\la_1,\la_2,\la_3}(M,M,M)=0.$
In the case of $k=1$, we have
\begin{eqnarray}\label{q=3-3}
0&=&(d\gamma)_{\la_1,\la_2,\la_3,\la}(L,M,M,L)|_{\la=0}\nonumber\\
&=&(\la_1-\la_2)\gamma_{0,\la_3,\la_1+\la_2}(L,M,M)-(\la_1-\la_3)\gamma_{0,\la_2,\la_1+\la_3}(L,M,M)\nonumber\\
&&-(\la_1+\la_2+\la_3)\gamma_{\la_1,\la_2,\la_3}(L,M,M).
\end{eqnarray}
Note that $\gamma_{\la_1,\la_2,\la_3}(L,M,M)$ is a homogeneous polynomial of degree $3$ and skew-symmetric in $\la_2$ and $\la_3$. Thus it is divisible by $\la_2-\la_3$. We can suppose that
\begin{eqnarray}\label{q=3-33}
\gamma_{\la_1,\la_2,\la_3}(L,M,M)
=(\la_2-\la_3)(a_1\la_1^2+a_2(\la_2^2+\la_3^2)+a_3\la_2\la_3
+a_4\la_1(\la_2+\la_3)),
\end{eqnarray}
where $a_1,a_2,a_3,a_4\in\CC$. Plugging \eqref{q=3-33} into \eqref{q=3-3} gives $a_4=0, a_3=2a_2,$ $a_1=-a_2$. Therefore,
\begin{eqnarray}
\gamma_{\la_1,\la_2,\la_3}(L,M,M)=a_2\phi_1,\ \mbox{where}\ \phi_1=(\la_2-\la_3)(\la_1+\la_2+\la_3)(-\la_1+\la_2+\la_3).
\end{eqnarray}
Note that $\phi_1$ is a coboundary of $\bar\g_{\la_1,\la_2}(M,M)=\la_2^2-\la_1^2$. In fact,
\begin{eqnarray*}
(d\bar\g)_{\la_1,\la_2,\la_3}(L,M,M)
&=&-(\la_1-\la_2)\g_{\la_1+\la_2,\la_3}(M,M)+(\la_1+\la_3)\g_{\la_1+\la_3,\la_2}(M,M)\\&=&
-(\la_1-\la_2)(\la_3^2-(\la_1+\la_2)^2)+(\la_1-\la_3)(\la_2^2-(\la_1+\la_3)^2)\\
&=&\phi_1.
\end{eqnarray*}
Similarly, suppose that
\begin{eqnarray}\label{q=3-5}
\gamma_{\la_1,\la_2,\la_3}(L,L,M)=(\la_1-\la_2)(b_1(\la_1^2+\la_2^2)+b_2\la_3^2+b_3(\la_1+\la_2)\la_3+b_4\la_1\la_2),
\end{eqnarray}
where $b_1,b_2,b_3,b_4\in\CC$. Substituting \eqref{q=3-5} into the following equality
\begin{eqnarray*}
0&=&(d\gamma)_{\la_1,\la_2,\la_3,\la}(L,L,M,L)|_{\la=0}\nonumber\\
&=&(\la_1-\la_2)\gamma_{\la_1+\la_2,0,\la_3}(L,L,M)+(\la_1-\la_3)\gamma_{\la_2,0,\la_1+\la_3}(L,L,M)\nonumber\\
&&-(\la_2-\la_3)\gamma_{\la_1,0,\la_2+\la_3}(L,L,M)-(\la_1+\la_2+\la_3)\gamma_{\la_1,\la_2,\la_3}(L,L,M)
\end{eqnarray*}
gives $b_4=b_1+b_2$. Hence,
\begin{eqnarray}\label{llm}
\gamma_{\la_1,\la_2,\la_3}(L,L,M)=(\la_1-\la_2)(b_1(\la_1^2+\la_2^2)+b_2\la_3^2+(b_1+b_2)\la_1\la_2+b_3(\la_1+\la_2)\la_3).
\end{eqnarray}
On the other hand, there is a 2-cochain $\bar\g_{\la_1,\la_2}(L,M)=b_1\la_1^2+b_2\la_1\la_2$ such that
\begin{eqnarray}\label{q=3-88}
(d\bar\g)_{\la_1,\la_2,\la_3}(L,L,M)+\gamma_{\la_1,\la_2,\la_3}(L,L,M)=-(b_1+b_2-b_3)(\la_1-\la_2)(\la_1+\la_2)\la_3.
\end{eqnarray}
Thus $\gamma_{\la_1,\la_2,\la_3}(L,L,M)$ in \eqref{llm} is equivalent to a constant factor of $\chi:=\chi_{\la_1,\la_2,\la_3}(L,L,M)=(\la_1-\la_2)(\la_1+\la_2)\la_3$, which is not a coboundary. By \cite[Theorem 7.1]{BKV},
$\Lambda_3:=\gamma_{\la_1,\la_2,\la_3}(L,L,L)=(\la_1-\la_2)(\la_1-\la_3)(\la_2-\la_3)$ (up to a constant factor) is a 3-cocycle, but not a coboundary. Therefore, ${\rm dim\,\t H}^3(\W,\CC)=2$. Specifically, ${\rm \t H}^3(\W,\CC)=\CC \chi \oplus \CC \Lambda_3.$

For $q=4$, three cases (i.e., $k=1, 2, 3$) should be taken into account. Let  $\g\in\t D^4(\W,\CC)$ be a $4$-cocycle. By using the method of undetermined coefficients and doing similar calculations to the case when $q=3$ , we obtain
 \begin{eqnarray}\label{q=4-1}
 \g_{\la_1,\la_2,\la_3,\la_4}(L,M,M,M)&=&c (\la_2-\la_3)(\la_3-\la_4)(\la_2-\la_4)(\la_2+\la_3+\la_4),\\
  \g_{\la_1,\la_2,\la_3,\la_4}(L,L,M,M)&=&(\la_1-\la_2)(\la_3-\la_4)(c_1(\la_1^2+\la_2^2)+c_2(\la_3+\la_4)^2\nonumber\\&&+(c_1+c_2)\la_1\la_2+c_3(\la_1+\la_2)(\la_3+\la_4)),\\
    \g_{\la_1,\la_2,\la_3,\la_4}(L,L,L,M)&=&(\la_1-\la_2)(\la_2-\la_3)(\la_1-\la_3)(e_1(\la_1+\la_2+\la_3)+e_2\la_4),
 \end{eqnarray}
where $c, c_1,c_2,c_3,e_1,e_2\in\CC$. And there exist three 3-cochains of degree 3
 \begin{eqnarray}
\bar\g_{\la_1,\la_2,\la_3}(M,M,M)&=&(\la_1-\la_2)(\la_1-\la_3)(\la_2-\la_3),\label{q=4-11} \\
\bar\g_{\la_1,\la_2,\la_3}(L,M,M)&=&(\la_2-\la_3)(c_1\la_1^2+c_2\la_1(\la_2+\la_3)),\\
\bar\g_{\la_1,\la_2,\la_3}(L,L,M)&=&(\la_1-\la_2)(\la_1^2+\la_2^2),
 \end{eqnarray}
such that
\begin{eqnarray}
\g_{\la_1,\la_2,\la_3,\la_4}(L,M,M,M)-c (d\bar\g)_{\la_1,\la_2,\la_3,\la_4}(L,M,M,M)&=&0,\label{q=4-2}\\
\gamma_{\la_1,\la_2,\la_3,\la_4}(L,L,M,M)+(d\bar\g)_{\la_1,\la_2,\la_3,\la_4}(L,L,M,M)&=&(c_3-c_1-c_2)\psi_1,\label{q=4-3}\\
\g_{\la_1,\la_2,\la_3,\la_4}(L,L,L,M)+e_2(d\bar\g)_{\la_1,\la_2,\la_3,\la_4}(L,L,L,M)&=&(e_1-e_2)\psi_2,\label{q=4-4}
 \end{eqnarray}
 where
 \begin{eqnarray}
&& \psi_1:=\psi_{1\,{\la_1,\la_2,\la_3,\la_4}}(L,L,M,M)=(\la_1-\la_2)(\la_1+\la_2)(\la_3-\la_4)(\la_3+\la_4),\label{q=4-5}\\ &&\psi_2:=\psi_{2\,{\la_1,\la_2,\la_3,\la_4}}(L,L,L,M)=(\la_1-\la_2)(\la_2-\la_3)(\la_1-\la_3)(\la_1+\la_2+\la_3).\label{q=4-6}
  \end{eqnarray}
 Moreover, $4\psi_1=-(d\bar\psi)_{\la_1,\la_2,\la_3,\la_4}(L,L,M,M)$ with $\bar\psi_{\la_1,\la_2,\la_3}(L,M,M)=(\la_2-\la_3)(3\la_1^2-(\la_2^2+\la_3^2))$.  This, together with \eqref{q=4-2}-- \eqref{q=4-6}, gives ${\rm \tilde H}^4(\W,\CC)=\CC \psi_2$.

For $q=5$, we need to consider $k=2$, $3$. Let $\g\in\t D^5(\W,\CC)$ be a $5$-cocycle. We obtain
\begin{eqnarray}
 \g_{\la_1,\la_2,\la_3,\la_4,\la_5}(L,L,M,M,M)&=&(\la_1-\la_2)(\la_3-\la_4)(\la_3-\la_5)(\la_4-\la_5)\nonumber\\
 && \times (\bar a_1(\la_1+\la_2)+\bar a_2(\la_3+\la_4+\la_5)),\label{q=5-1}\\
 \g_{\la_1,\la_2,\la_3,\la_4,\la_5}(L,L,L,M,M)&=&(\la_1-\la_2)(\la_1-\la_3)(\la_2-\la_3)(\la_4-\la_5)\nonumber\\
 && \times (\bar b_1(\la_1+\la_2+\la_3)+\bar b_2(\la_4+\la_5)),\label{q=5-2}
\end{eqnarray}
where $\bar a_1,\bar a_2,\bar b_1,\bar b_2\in\CC.$ On the other hand, there exist two $4$-cochains of degree 4
\begin{eqnarray*}
\bar\g_{\la_1,\la_2,\la_3,\la_4}(L,M,M,M)&=&\la_1(\la_2-\la_3)(\la_3-\la_4)(\la_2-\la_4),\\
\bar\g_{\la_1,\la_2,\la_3,\la_4}(L,L,M,M)&=&(\la_1-\la_2)(\la_3-\la_4)(\la_1^2+\la_2^2),
\end{eqnarray*}
such that
\begin{eqnarray}
\g_{\la_1,\la_2,\la_3,\la_4,\la_5}(L,L,M,M,M)+a_1(d\bar\g)_{\la_1,\la_2,\la_3,\la_4,\la_5}(L,L,M,M,M)&=& (a_2-a_1)\varphi_1,\label{q=5-11}\\
\g_{\la_1,\la_2,\la_3,\la_4,\la_5}(L,L,L,M,M)+b_1(d\bar\g)_{\la_1,\la_2,\la_3,\la_4,\la_5}(L,L,L,M,M)
 &=&(b_2-b_1)\varphi_2,\label{q=5-22}
\end{eqnarray}
where \begin{eqnarray}\label{q=5-22-1}
\varphi_1:&=&\varphi_{1\,{\la_1,\la_2,\la_3,\la_4,\la_5}}(L,L,M,M,M)\nonumber\\ &=&(\la_1-\la_2)(\la_3-\la_4)(\la_3-\la_5)(\la_4-\la_5)(\la_3+\la_4+\la_5),\\ \varphi_2:&=&\varphi_{2\,{\la_1,\la_2,\la_3,\la_4,\la_5}}(L,L,L,M,M)\nonumber\\ &=&(\la_1-\la_2)(\la_1-\la_3)(\la_2-\la_3)(\la_4-\la_5)(\la_4+\la_5).
\end{eqnarray}
Furthermore, there exists another one $4$-cochains of degree 4
\begin{eqnarray*}
\bar\varphi_{\la_1,\la_2,\la_3,\la_4}(L,L,M,M)=(\la_1-\la_2)(\la_3-\la_4)(\la_1\la_2-\la_3\la_4),
\end{eqnarray*}
such that $2\varphi_2=(d\bar\varphi)_{\la_1,\la_2,\la_3,\la_4,\la_5}(L,L,L,M,M)$, namely, $\varphi_2$ is a coboundary. By \eqref{q=5-11} and \eqref{q=5-22}, ${\rm dim\, \t H}^5(\W,\CC)=1$, and ${\rm \tilde H}^5(\W,\CC)=\CC \varphi_1$.

For $q=6$, it only needs to consider the case when $k=3$. One can check that
\begin{eqnarray}\label{q=6}
\Lambda:&=&\g_{\la_1,\la_2,\la_3,\la_4,\la_5,\la_6}(L,L,L,M,M,M)\nonumber\\&=&(\la_1-\la_2)(\la_2-\la_3)(\la_1-\la_3)(\la_4-\la_5)(\la_4-\la_6)(\la_5-\la_6)
\end{eqnarray}
is a 6-cocycle. It is not a coboundary. Because it can be the coboundary of a 5-cochain of degree 5, which must be a constant factor of  $ \g_{\la_1,\la_2,\la_3,\la_4,\la_5}(L,L,M,M,M)$ in \eqref{q=5-1}, whose coboundary is zero. Therefore, ${\rm dim\,\t H}^6(\W,\CC)=1$ and ${\rm \t H}^6(\W,\CC)=\CC \Lambda$. This proves \eqref{main1}.

 According to \cite[Proposition 2.1]{BKV}, the map $\gamma\mapsto\partial \gamma$ gives an isomorphism ${\rm \t H}^q(\W,\CC)\cong {\rm H}^q(\pa\t C^\bullet) $ for $q\geq 1$.
Therefore
\begin{eqnarray}\label{pa-c}
{\rm H}^q(\pa\t C^\bullet)=\left\{
\begin{array}{ll}
\CC (\pa\chi) \oplus \CC (\pa\Lambda_3) &{if}\ q=3,\\
\CC (\pa\psi_2) &{if}\ q=4,\\
\CC (\pa\varphi_1) &{if}\ q=5,\\
\CC (\pa\Lambda) &{if}\ q=6,\\
0 &{otherwise}.
\end{array}
\right.
\end{eqnarray}

It remains to compute ${\rm H}^\bullet(\W,\CC)$. This is based on
the short exact sequence of complexes
 \begin{equation}\label{exact}
\begin{CD}
0@>>> {\pa\t C^\bullet} @>{\rm \iota}>> {\t C^\bullet} @>{\rm \pi}>> C^\bullet @>>> 0
\end{CD}
\end{equation}
 where $\iota$  and $\rm \pi$ are the embedding and the natural projection, respectively. The exact sequence \eqref{exact} gives the following long exact sequence of cohomology groups (cf. \cite{BKV}):
\begin{equation}\label{key}
\begin{CD}
\cdots @>>> {\rm H}^q(\pa\t C^\bullet) @>{\rm \iota}_q>> {\rm \t H}^q(\W,\CC)  @>{\rm \pi}_q>> {\rm H}^q(\W,\CC) @>{\rm \omega}_q>>\\
@>>> {\rm H}^{q+1}(\pa\t C^\bullet) @>{\rm \iota}_{q+1}>> {\rm \t H}^{q+1}(\W,\CC)  @>{\rm \pi}_{q+1}>> {\rm H}^{q+1}(\W,\CC) @>>>\cdots \\
\end{CD}
\end{equation}
where $\iota_q, \pi_{q}$ are induced by $\iota, \pi$ respectively and $w_q$ is the $q-$th connecting homomorphism.
Given $\partial \gamma\in {\rm H}^q(\pa\t C^\bullet)$ with a nonzero element $\gamma\in {\rm \t H}^q(\W,\CC)$,
then $\iota_q(\partial \gamma)=\partial \gamma\in {\rm \t H}^q(\W,\CC)$.
Since ${\rm deg\,}(\partial\gamma)={\rm deg\,}(\gamma)+1=q+1$,  we have $\partial\gamma=0\in {\rm \t H}^q(\W,\CC)$.
Then the image of $\iota_q$ is zero for any $q\in\Z_+$.
Because ${\rm ker}({\rm \pi}_{q})={\rm im}({\rm \iota}_{q})=\{0\}$ and ${\rm im}({\rm \omega}_{q})={\rm ker}({\rm \iota}_{q+1})={\rm H}^{q+1}(\pa\t C^\bullet)$,
we obtain the following short exact sequence
\begin{equation}\label{key1}
\begin{CD}
0 @>>> {\rm \t H}^q(\W,\CC)  @>{\rm \pi}_{q}>> {\rm H}^q(\W,\CC) @>{\rm \omega}_q>>{\rm H}^{q+1}(\pa\t C^\bullet)  @>>>0.
\end{CD}
\end{equation}
Therefore,
\begin{eqnarray}\label{key111}
{\rm dim\, H}^q(\W,\CC)={\rm dim\,\t H}^q(\W,\CC)+{\rm dim\, H}^{q+1}(\pa\t C^\bullet), \ \mbox{for\ all} \ q\geq 0.
\end{eqnarray}
Then \eqref{main2} follows from \eqref{key111}.
Moreover, we can give a basis for ${\rm H}^q(\W,\CC)$. Indeed,
any basis of ${\rm H}^q(\W,\CC)$ can be obtained by combining the images of a basis of ${\rm \t H}^q(\W,\CC)$
with the pre-images of a basis of ${\rm \t H}^{q+1}(\W,\CC)$. For a nonzero $\partial\varphi\in{\rm H}^{q+1}(\pa\t C^\bullet)$ with $\varphi$ being a $(q+1)$-cocycle, \eqref{7++2} gives
\begin{eqnarray}
d(\tau(\pa\varphi))=(d\tau+\tau d)(\pa\varphi)=({\rm deg\,} (\pa\varphi)-(q+1))(\pa\varphi)=((q+2)-(q+1))(\pa\varphi)=\pa\varphi.
\end{eqnarray}
Thus the pre-image of $\pa\varphi$ under the connecting homomorphism $\omega_p$ is $\omega_q^{-1}(\partial\varphi)=\tau(\pa\varphi)$.

Finally, we finish our proof by giving a basis of ${\rm H}^q(\W,\CC)$ for $q=2,\dots,6$.
For $q=2$, we have known that ${\rm \t H}^2(\W,\CC)=0$ and ${\rm H}^3(\pa\t C^\bullet)=\CC (\pa\chi) \oplus \CC (\pa\Lambda_3)$. By \eqref{7++} and \eqref{q=3-88},
\begin{eqnarray*}
\bar\chi:&=&(\tau(\pa\chi))_{\la_1,\la_2}(L,M)\\&=&(-1)^2\frac{\pa}{\pa\la}(\pa\chi)_{\la_1,\la_2,\la}(L,M,L)|_{\la=0}\\
&=&-\frac{\pa}{\pa\la}(\la_1+\la_2+\la)(\la_1^2-\la^2)\la_2|_{\la=0}\\
&=&-\la_1^2\la_2,\\
\bar\Lambda_3:&=&(\tau(\pa\Lambda_3))_{\la_1,\la_2}(L,L)\\&=&(-1)^2\frac{\pa}{\pa\la}(\pa\Lambda_3)_{\la_1,\la_2,\la}(L,L,L)|_{\la=0}\\
&=&\frac{\pa}{\pa\la}(\la_1+\la_2+\la)(\la_1-\la_2)(\la_2-\la)(\la_1-\la)|_{\la=0}\\
&=&-\la_1^3+\la_2^3.
\end{eqnarray*}
This gives ${\rm H}^2(\W,\CC)=\CC \bar \chi \oplus \CC \bar \Lambda_3.$
For $q=3$, by \eqref{7++}, \eqref{q=4-6} and \eqref{pa-c},
\begin{eqnarray*}
\bar \psi:&=&(\tau(\pa\psi_2))_{\la_1,\la_2,\la_3}(L,L,M)\\&=&(-1)^3\frac{\pa}{\pa\la}(\pa\psi_2)_{\la_1,\la_2,\la_3,\la}(L,L,M,L)|_{\la=0}\\
&=&\frac{\pa}{\pa\la}(\la_1+\la_2+\la_3+\la)(\la_1-\la_2)(\la_1-\la)(\la_2-\la)(\la_1+\la_2+\la)|_{\la=0}\\
&=&-\la_1^4-\la_1^3\la_3+\la_2^3(\la_2+\la_3).
\end{eqnarray*}
Hence, ${\rm H}^3(\W,\CC)={\rm \t H}^3(\W,\CC)\oplus \CC \bar\psi=\CC \chi \oplus \CC \Lambda_3 \oplus \CC \bar\psi$. By \eqref{7++}, \eqref{q=5-22-1}, \eqref{q=6} and \eqref{pa-c},
\begin{eqnarray*}
\bar \varphi:&=&(\tau(\pa\varphi_1))_{\la_1,\la_2,\la_3,\la_4}(L,M,M,M)\\
&=&(-1)^4\frac{\pa}{\pa\la}(\pa\varphi_1)_{\la_1,\la_2,\la_3,\la_4,\la}(L,M,M,M,L)|_{\la=0}\\
&=&-\frac{\pa}{\pa\la}(\la_1+\la_2+\la_3+\la_4+\la)(\la_1-\la)(\la_3-\la_2)(\la_3-\la_4)(\la_4-\la_2)(\la_2+\la_3+\la_4)|_{\la=0}\\
&=&(\la_2-\la_3)(\la_2-\la_4)(\la_3-\la_4)(\la_2+\la_3+\la_4)^2,\\
\bar\Lambda:&=&(\tau(\pa\Lambda))_{\la_1,\la_2,\la_3,\la_4,\la_5}(L,L,M,M,M)\\
&=&(-1)^5\frac{\pa}{\pa\la}(\pa\Lambda)_{\la_1,\la_2,\la_3,\la_4,\la_5,\la}(L,L,M,M,M,L)|_{\la=0}\\
&=&\frac{\pa}{\pa\la}\big(\mbox{$\sum_{i=1}^5$}\la_i+\la\big)(\la_1-\la_2)(\la_1-\la)(\la_2-\la)(\la_4-\la_5)(\la_4-\la_3)(\la_5-\la_3)|_{\la=0}\\
&=&-(\la_1-\la_2)(\la_3-\la_4)(\la_3-\la_5)(\la_4-\la_5)\big(\la_1^2+\la_1\la_2+\la_2^2+(\la_1+\la_2)(\la_3+\la_4+\la_5)\big)\\
&\equiv&\la_1\la_2(\la_1-\la_2)(\la_3-\la_4)(\la_3-\la_5)(\la_4-\la_5)(\mbox{mod}\
\pa\t C^5(\W, \CC)).
\end{eqnarray*}
Therefore, ${\rm H}^4(\W,\CC)=\CC \psi_2\oplus \CC \bar\varphi$,  ${\rm H}^5(\W,\CC)=\CC \varphi_1\oplus \CC \bar\Lambda $ and
 ${\rm H}^6(\W,\CC)=\CC \Lambda.$ Thus (1) is proved.

(2) Define an operator $\tau_2:\t C^q(\W,\CC_a)\rightarrow
\t C^{q-1}(\W,\CC_a)$ by
\begin{eqnarray}\label{7-3}
(\tau_2
\g)_{\la_1,\cdots,\la_{q-1}}(X_1,\cdots,X_{q-1})=(-1)^{q-1}\g_{\la_1,\cdots,\la_{q-1},\la}(X_1,\cdots,X_{q-1},L)|_{\la=0},
\end{eqnarray}
for $X_1,\cdots,X_{q-1}\in\{L,M\}$. By the fact that $\pa\t C^q(\W,\CC_a)=(a+\sum_{i=1}^q\la_i)\t
C^q(\W,\CC_a)$, we
have
\begin{eqnarray}\label{7-4}
((d\tau_2+\tau_2 d)
\g)_{\la_1,\cdots,\la_{q}}(X_1,\cdots,X_q)&=&\big(\mbox{$\sum_{i=1}^q$}\la_i\big)\g_{\la_1,\cdots,\la_{q}}(X_1,\cdots,X_q)\nonumber\\
&\equiv& -a \g_{\la_1,\cdots,\la_{q}}(X_1,\cdots,X_q) \ (\mbox{mod}\
\pa\t C^q(\W, \CC_a)).
\end{eqnarray}
Let $\g\in\t C^q(\mathcal{W},\CC_a)$ be
a $q$-cochain such that $d\g\in \pa\t C^{q+1}(\W,\CC_a)$, namely,
there is a $(q+1)$-cochain $\phi$ such that
$d\g=(a+\sum_{i=1}^{q+1}\la_i)\phi$. By \eqref{7-3}, $\tau_2 d\g=(a+\sum_{i=1}^q\la_i)\tau_2\phi\in\pa\t
C^q(\W,\CC_a)$. It follows from \eqref{7-4} that
$\g\equiv-d(a^{-1}\tau_2\g)$ is a reduced coboundary. This proves (2).

(3) In this case, $\pa\t
C^q(\W,M_{\Delta,\a})=(\pa+\sum_{i=1}^q\la_i)\t C^q(\W,
M_{\Delta,\a})$. As in the proof of (2), we define an
operator $\tau_3:C^q(\W,M_{\Delta,\a})\rightarrow
C^{q-1}(\W,M_{\Delta,\a})$ by
\begin{eqnarray*}\label{5-3}
(\tau_3
\g)_{\la_1,\cdots,\la_{q-1}}(X_1,\cdots,X_{q-1})=(-1)^{q-1}\g_{\la_1,\cdots,\la_{q-1},\la}(X_1,\cdots,X_{q-1},L)|_{\la=0},
\end{eqnarray*}
for $X_1,\cdots,X_{q-1}\in\{L,M\}$. Then
\begin{eqnarray}\label{5-4}
((d\tau_3+\tau_3 d)
\g)_{\la_1,\cdots,\la_{q}}(X_1,\cdots,X_q)&=&L_\la\g_{\la_1,\cdots,\la_{q}}(X_1,\cdots,X_{q})|_{\la=0}+
\big(\mbox{$\sum\limits_{i=1}^q$}\la_i\big)\g_{\la_1,\cdots,\la_{q}}(X_1,\cdots,X_q)\nonumber\\[-6pt]
&=&\big(\pa+\a+\mbox{$\sum\limits_{i=1}^q$}\la_i\big)\g_{\la_1,\cdots,\la_{q}}(X_1,\cdots,X_q)\nonumber\\[-6pt]
&\equiv& \a \g_{\la_1,\cdots,\la_{q}}(X_1,\cdots,X_q) \ (\mbox{mod}\
\pa\t C^q(\W, M_{\Delta,\a})).
\end{eqnarray}
If $\g$ is a reduced $q$-cocycle, it follows from \eqref{5-4}
that $\g\equiv d(\a^{-1}\tau_3\g)$ is a reduced coboundary, since $\a\neq
0$. Thus $ {\rm H}^q(\mathcal{W},M_{\Delta,\a})=0$ for all $q\geq 0$.

This completes the proof of Theorem \ref{thm-5}. \QED
\begin{rema}\rm  Denote by ${\rm Lie}(\mathcal{W})_-$ the annihilation Lie algebra of $\mathcal{W}$.
Note that ${\rm Lie}(\mathcal{W})_-$ is isomorphic to the subalgebra spanned by
$\{L_n, M_n\big| -1\leq n\in\Z\}$ of the centerless
$W$-algebra $W(2,2)$. By \cite[Corollary 6.1]{BKV}, ${\rm \t H}^q(\mathcal{W},\CC)\cong{\rm H}^q({\rm Lie}(\mathcal{W})_-,\CC)$. Thus we have determined
the cohomology of ${\rm Lie}(\mathcal{W})_-$ with trivial coefficients.
\end{rema}

\vspace{4mm} \noindent\bf{\footnotesize Acknowledgements.}\ \rm
{\footnotesize The authors would like to thank
the referees for helpful suggestions.
This work was supported by National Natural Science
Foundation grants of China (11301109, 11526125) and the Research Fund for the Doctoral Program of Higher Education (20132302120042).}\\
\vskip18pt \small\footnotesize
\parskip0pt\lineskip1pt


\parskip=0pt\baselineskip=1pt
\begin{thebibliography}{9999}
\def\RE#1{\bibitem{#1}\label{#1}}

\bibitem {K1} Kac V.,  Vertex algebras for beginners, Univ. Lect. Series 10, AMS,
1996

\bibitem {K3} Kac V., The idea of locality, in: H.-D. Doebner, et al. (Eds.),
Physical Applications and Mathematical Aspects of Geometry, Groups
and Algebras, World Sci. Publ., Singapore, 1997, 16--32


\bibitem {DK} D'Andrea A., Kac V., Structure theory of finite
conformal algebras, Sel. Math., New Ser., 1998, 4, 377--418

\bibitem {FK} Fattori D., Kac V., Classification of finite simple Lie
conformal superalgebras, J. Algebra, 2002, 258(1), 23--59

\bibitem {CK} Cheng S.-J., Kac V., Conformal modules, Asian
J. Math., 1997, 1(1), 181--193

\bibitem{BKV} Bakalov B., Kac V., Voronov A., Cohomology of
conformal algebras, Comm. Math. Phys., 1999, 200, 561--598

\bibitem {FSW} Fan G., Su Y., Wu H., Loop Heisenberg-Virasoro Lie conformal algebra, J. Math. Phys., 2014, 55, 123508


\bibitem{SY} Su Y., Yuan L., Schr\"{o}dingger-Virasoro Lie conformal
algebra, J. Math. Phys., 2013, 54, 053503

\bibitem{WCY} Wu H., Chen Q., Yue X., Loop Virasoro Lie conformal algebra, J. Math. Phys., 2014, 55, 011706

\bibitem{ZD} Zhang W., Dong C., $W$-algebra $W(2, 2)$ and the vertex
operator algebra $L(\frac1 2, 0)\otimes L(\frac1 2, 0)$, Comm.
Math. Phys., 2009, 285, 991--1004


\bibitem{XY} Xu Y., Yue X., $W(a,b)$ Lie conformal algebra and its conformal
module of rank one, Algebra Colloq., 2015, 22(3), 405--412


\bibitem{MP} Mikhalev A.V., Pinchuk I.A., Universal Central extensions of Lie conformal algebras,
Part 2: supercase, Moscow University Mathematics Bulletin, 2010, 65(1), 34--38

\end{thebibliography}
\end{document}